\documentclass[11pt]{article}
\usepackage{amsmath}
\usepackage{amsfonts}
\usepackage{amssymb}
\usepackage{mathrsfs}
\usepackage{latexsym}
\usepackage[dvips]{graphicx}
\usepackage{multirow}
\newtheorem{theorem}{\bf Theorem}[section]
\newtheorem{lem}[theorem]{\bf Lemma}
\newtheorem{rem}[theorem]{\bf Remark}

\def \bG {\Bbb G}

\begin{document}

\title{{\bf Sharp bounds for Hardy type operators on higher-dimensional product spaces}}
\author{{\bf QianJun HE  \quad \bf DunYan YAN\footnote{{Corresponding author}} }
\\\footnotesize\textit{School of Mathematics, Graduate University, Chinese Academy of Sciences, Beijing 100049, China}
\\\footnotesize\textit{ heqianjun16@mails.ucas.ac.cn \quad\quad ydunyan@ucas.ac.cn
}}

\date{}
\maketitle \vspace{0.2cm}

\noindent{{\bf Abstract}} \quad In this paper, we investigate  a class of fractional Hardy type operators $\mathscr{H}_{\beta_{1},\cdots,\beta_{m}}$ defined on  higher-dimensional product spaces $\mathbb{R}^{n_{1}}\times\mathbb{R}^{n_{2}}\times\cdots\times\mathbb{R}^{n_{m}}$. We use novel methods to obtain two main results. One is that the operator $\mathscr{H}_{\beta_{1},\cdots,\beta_{m}}$ is bounded from $L^{p}(\mathbb{R}^{n_{1}}\times\mathbb{R}^{n_{2}}\times\cdots\times\mathbb{R}^{n_{m}},|x|^{\gamma})$ to $L^{q}(\mathbb{R}^{n_{1}}\times\mathbb{R}^{n_{2}}\times\cdots\times\mathbb{R}^{n_{m}},|x|^{\alpha})$ and the sharp bound of the operator $\mathscr{H}_{\beta_{1},\cdots,\beta_{m}}$ is worked out. The other is that when $\alpha=\gamma=(0,\cdots,0)$, the norm of the operator $\mathscr{H}_{\beta_{1},\cdots,\beta_{m}}$ is obtained.
\medskip

\noindent{{\bf 2010 Subject Classification}} \quad 42B20, 42B25


\section{\hspace{-0.3cm}{\bf} Introduction }
Let $f$ be non-negative measurable function on the $m$-fold product space $\mathbb{R}^{n_{1}}\times\mathbb{R}^{n_{2}}\times\cdots\times\mathbb{R}^{n_{m}}$, the Hardy type operator $\mathscr{H}_{\beta_{1},\beta_{2},\cdots,\beta_{m}}$ is defined by
\begin{equation}\label{eq:1}
\mathscr{H}_{\beta_{1},\cdots,\beta_{m}}(f)(x):=\prod_{i=1}^{m}\frac{1}{|B(0,|x_{i}|)|^{1-\frac{\beta_{i}}{n_{i}}}}\int_{|y_{1}|
<|x_{1}|}\cdots\int_{|y_{m}|<|x_{m}|}f(y_{1},\cdots,y_{m})dy_{1}\cdots dy_{m},
\end{equation}
where $x=(x_    {1},x_{2},\cdots,x_{m})\in\mathbb{R}^{n_{1}}\times\mathbb{R}^{n_{2}}\times\cdots\times\mathbb{R}^{n_{m}}$ with $\prod_{i=1}^{m}|x_{i}|\neq0$ and $0\leq\beta_{i}<n_{i}$ for $i=1,2,\cdots,m$.

Obviously, the  operator $\mathscr{H}_{\beta_{1},\beta_{2},\cdots,\beta_{m}}$ is natural generalization of
the classical Hardy operators, such as the operator $H$ \cite{H1,H2}, the fractional Hardy operator $\mathbb{H}_{\beta}$ \cite{LZ}, the Hardy operator $\mathcal{H}$ on the $n$-dimensional product space \cite{WLY} and the Hardy type operator $\mathcal{H}_{m}$ on $m$-dimensional product spaces \cite{LYZ}.

 If $f$ is a non-negative measurable function on $\mathbb{G}=(0,\infty)$, the classical Hardy operator defined as
$$H(f)(x):=\frac{1}{x}\int_{0}^{x}f(t)dt,\quad x>0.$$
The following Theorem A due to Hardy \cite{H1,H2} is well known.

{\noindent}\textrm{\bf Theorem A} \quad If $f$ is a non-negative measurable function on $\mathbb{G}$, let $1<p<\infty$ and $\alpha<p-1$, then the following two inequalities
	$$\|H(f)\|_{L^{p}(\mathbb{G})}\leq\frac{p}{p-1}\|f\|_{L^{p}(\mathbb{G})} \quad\text{and}\quad \|H(f)\|_{L^{p}(\mathbb{G},x^\alpha)}\leq\frac{p}{p-1-\alpha}\|f\|_{L^{p}(\mathbb{G},x^{\alpha})}$$
hold, and both constants $\frac{p}{p-1}$ and $\frac{p}{p-1-\alpha}$ are sharp.

Recall that, for a nonnegative measurable function $f$ on $\mathbb{R}^{n}$, the $n$-dimensional fractional Hardy operator $\mathbb{H}_{\beta}$ with spherical mean is defined by
\begin{equation}\label{eq:2}
\mathbb{H}_{\beta}(f)(x)=\frac{1}{|B(0,|x|)|^{1-\frac{\beta}{n}}}\int_{|y|<|x|}f(y)dy,
\end{equation}
where $x\in\mathbb{R}^{n}\backslash\{0\}$ and $0\leq\beta<n$.

In 2015, Lu and Zhao considered the operator defined by $\eqref{eq:2}$ and obtained the following
Theorem B.

{\noindent}\textrm{\bf Theorem B} \quad If $f$ is a nonnegative measurable function on $\mathbb{R}^{n}$, let  $0<\beta<n$, $1<p<q<\infty$ and $\frac{1}{q}=\frac{1}{p}-\frac{\beta}{n}$, then the inequality
\begin{equation}\label{eq:3}
\|\mathbb{H}_{\beta}(f)\|_{L^{q}}\leq\left(\frac{p^{\prime}}{q}\right)^{1/q}\left(\frac{n}{q\beta}\cdot B\left(\frac{n}{q\beta},\frac{n}{q^{\prime}\beta}\right)\right)^{-\beta/n}\|f\|_{L^{p}}
\end{equation}
holds, where the constant of inequality in $\eqref{eq:3}$ is sharp.

For the weighted case, in 1985, Sawyer \cite{S} considered the weighted Hardy inequalities with general
weight functions $u$ and $v$ only on two-dimensional product space. However, it is much hard to apply the method in \cite{S} to the case dimensional greater than two. In \cite{WLY}, Wang, Lu and Yan studied the Hardy operator $\mathcal{H}$ with power weight on the $m$-dimensional case, and obtained the following result.

{\noindent}\textrm{\bf Theorem C} \quad If $f$ is a nonnegative measurable function on $\bG^m$ and $x=(x_{1},x_{2},\cdots,x_{m})\in\bG^m$, let $\alpha=(\alpha_{1},\alpha_{2},\cdots,\alpha_{m})$, $\beta=(\beta_{1},\beta_{2},\cdots,\beta_{m})$, $x^{\alpha}=x_{1}^{\alpha_{1}}x_{2}^{\alpha_{2}}\cdots x_{m}^{\alpha_{m}}$, $x^{\beta}=x_{1}^{\beta_{1}}x_{2}^{\beta_{2}}\cdots x_{m}^{\beta_{m}}$, $1<p\leq q<\infty$, $\frac{1}{q}+1=\frac{1}{p}+\frac{1}{r}$, $\alpha_{i}<p-1$ and $\frac{\beta_{i}+1}{q}=\frac{\alpha_{i}+1}{p}$ with $1\leq i\leq m$, then
\begin{equation}\label{same power}
\|\mathcal{H}f\|_{L^{p}(\bG^m,x^{\alpha})}\leq \left(\prod_{i=1}^{m}\frac{p}{p-\alpha_{i}-1}\right)\|f\|_{L^{p}(\bG^m,x^{\alpha})}
\end{equation}
and
\begin{equation}\label{different power}
\|\mathcal{H}f\|_{L^{q}(\bG^m,x^{\beta})}\leq \left(\prod_{i=1}^{m}\frac{q}{r(q-\beta_{i}-1)}\right)^{\frac{1}{r}}\|f\|_{L^{p}(\bG^m,x^{\alpha})},
\end{equation}
 where the operator $\mathcal{H}$ is the fractional Hardy type operator $\mathscr{H}_{\beta_{1},\beta_{2},\cdots,\beta_{m}}$ with $\beta_{1}=\cdots=\beta_{m}=0$ and $n_{1}=\cdots=n_{m}=1$.
The authors \cite{WLY} not sure the constant $\left(\prod_{i=1}^{m}\frac{q}{r(q-\beta_{i}-1)}\right)^{\frac{1}{r}}$ is sharp. However, it needs to be proved. Our results in the following will show that sharp constant
in $\eqref{different power}$ is not equal to $\left(\prod_{i=1}^{m}\frac{q}{r(q-\beta_{i}-1)}\right)^{\frac{1}{r}}$.

For Another high dimensional case with power weight, Lu, Yan and Zhao in \cite{LYZ} studied the operator $\mathcal{H}_{m}$ is definition of $\mathscr{H}_{\beta_{1},\beta_{2},\cdots,\beta_{m}}$ with $\beta_{1}=\beta_{2}=\cdots=\beta_{m}=0$ and obtained sharp constant $\prod_{j=1}^{m}\frac{p}{p-1-\alpha_{j}/n_{j}}$.

A natural question is to consider the case of higher-dimensional product space for fractional Hardy type operator $\mathscr{H}_{\beta_{1},\beta_{2},\cdots,\beta_{m}}$. We will use novel methods and ideas to study Hardy operator with power weight on higher-dimensional product spaces. For more information about the Hardy type operator, we refer to (\cite{B},\cite{BPN},\cite{KP},\cite{M2},\cite{PS}) and references therein.

Now, we formulate our main results as follows.
\begin{theorem}\label{main:1}
	Suppose that $1<p<q<\infty$, $0<\beta_{i}<n_{i}$ and $\frac{1}{q}=\frac{1}{p}-\frac{\beta_{i}}{n_{i}}$ with $1\leq i\leq m$. If $f\in L^{p}(\mathbb{R}^{n_{1}}\times\mathbb{R}^{n_{2}}\times\cdots\times\mathbb{R}^{n_{m}})$, then we have
	\begin{equation}\label{eq:4}
	\|\mathscr{H}_{\beta_{1},\beta_{2},\cdots,\beta_{m}}(f)\|_{L^{q}(\mathbb{R}^{n_{1}}\times\mathbb{R}^{n_{2}}\times\cdots\times\mathbb{R}^{n_{m}})}\leq C\|f\|_{L^{p}(\mathbb{R}^{n_{1}}\times\mathbb{R}^{n_{2}}\times\cdots\times\mathbb{R}^{n_{m}})}.
	\end{equation}
	Moreover,
	$$\|\mathscr{H}_{\beta_{1},\beta_{2},\cdots,\beta_{m}}\|_{L^{p}(\mathbb{R}^{n_{1}}\times\mathbb{R}^{n_{2}}\times\cdots\times\mathbb{R}^{n_{m}})\rightarrow L^{q}(\mathbb{R}^{n_{1}}\times\mathbb{R}^{n_{2}}\times\cdots\times\mathbb{R}^{n_{m}})}=C,$$
	where
	$$C=\prod_{i=1}^{m}\left(\frac{p^{\prime}}{q}\right)^{1/q}\left(\frac{n_{i}}{q\beta_{i}}\cdot B\left(\frac{n_{i}}{q\beta_{i}},\frac{n_{i}}{q^{\prime}\beta_{i}}\right)\right)^{-\beta_{i}/n_{i}}.$$
\end{theorem}

For two differences power weight, we have following result.
\begin{theorem}\label{main:2}
Suppose that $1<p\leq q<\infty$, $m\in\mathbb{N}$, $n_{i}\in\mathbb{N}$, $\gamma=(\gamma_{1},\gamma_{2},\cdots,\gamma_{m})$, $\alpha=(\alpha_{1},\alpha_{2},\cdots,\alpha_{m})$, $x_{i}\in\mathbb{R}^{n_{i}}$, $0\leq\beta_{i}<n_{i}$, $\gamma_{i}<n_{i}(p-1)$ and $\beta_{i}+\frac{\alpha_{i}+n_{i}}{q}=\frac{\gamma_{i}+n_{i}}{p}$ with $1\leq i\leq m$. If $f\in L^{p}(\mathbb{R}^{n_{1}}\times\mathbb{R}^{n_{2}}\times\cdots\times\mathbb{R}^{n_{m}},|x|^{\gamma})$ with $|x|^{\gamma}=|x_{1}|^{\gamma_{1}}|x_{2}|^{\gamma_{2}}\cdots|x_{m}|^{\gamma_{m}}$, then we have
	\begin{equation}\label{eq:5}
	\|\mathscr{H}_{\beta_{1},\beta_{2},\cdots,\beta_{m}}(f)\|_{L^{q}(\mathbb{R}^{n_{1}}\times\mathbb{R}^{n_{2}}
\times\cdots\times\mathbb{R}^{n_{m}},|x|^{\alpha})}\leq C^{*}\|f\|_{L^{p}(\mathbb{R}^{n_{1}}\times\mathbb{R}^{n_{2}}\times\cdots\times\mathbb{R}^{n_{m}},|x|^{\gamma})},
	\end{equation}
		Moreover,
		$$\|\mathscr{H}_{\beta_{1},\beta_{2},\cdots,\beta_{m}}\|_{L^{p}(\mathbb{R}^{n_{1}}\times\mathbb{R}^{n_{2}}
\times\cdots\times\mathbb{R}^{n_{m}},|x|^{\alpha})\rightarrow L^{q}(\mathbb{R}^{n_{1}}\times\mathbb{R}^{n_{2}}\times\cdots\times\mathbb{R}^{n_{m}},|x|^{\gamma})}=C^{*},$$
		where $C^{*}$ is equal to
		$$\prod_{i=1}^{m}\left|S^{n_{i}-1}\right|^{\frac{1}{q}-\frac{1}{p}+\frac{\beta_{i}}{n_{i}}}n_{i}^{\frac{1}{p}-\frac{1}{q}-\frac{\beta_{i}}{n_{i}}}
\left(\frac{n_{i}(p-1)}{n_{i}(p-1)-\gamma_{i}}\right)^{\frac{1}{p^{\prime}}+\frac{1}{q}}\left(\frac{p^{\prime}}{q}\right)^{1/q}\left(\frac{p}{q-p}\cdot B\left(\frac{p}{q-p},\frac{pq}{q^{\prime}(q-p)}\right)\right)^{\frac{1}{q}-\frac{1}{p}}.$$
\end{theorem}

It is worth mentioning that that proof in \cite{WLY} is not suitable to the  operator $\mathscr{H}_{\beta_{1},\beta_{2},\cdots,\beta_{m}}$. Although the idea in the paper is motivated  by the reference \cite{PS}, there are some essential differences. The difficulty is how to deal with the product space case. In this paper we will use the novel method to become as result in \cite{LZ}. The reconstruct some auxiliary functions to achieve the sharp bounds, which is quite different from \cite{PS}.

Throughout the note, we use the following notation. The definition of the usual beta function is defined by
$$B(z,w)=\int_{0}^{0}t^{z-1}(1-t)^{w-1}dt,$$ where $z$ and $w$ are complex numbers with the positive real parts. The set $B(0,|x|)$ denotes a open ball with center at the original point and radius $|x|$, and $|B(0,|x|)|$ denotes the volume of the ball $B(0,|x|)$. For one $m$-dimensional vector $\alpha=(\alpha_{1},\alpha_{2},\cdots,\alpha_{m})$ and $|x|^{\alpha}=|x_{1}|^{\alpha_{1}}|x_{2}|^{\alpha_{2}}\cdots|x_{m}|^{\alpha_{m}}$, $x=(x_{1},x_{2},\cdots,x_{m})\in\mathbb{R}^{n_{1}}\times\mathbb{R}^{n_{2}}\times\cdots\times\mathbb{R}^{n_{m}}$.  For a real number $p$, $1<p<\infty$, $p^{\prime}$ is the conjugate number of $p$, that is, $1/p+1/p^{\prime}=1$.

\section{Preliminaries}

\qquad To reduce the dimension of function space, we need the following lemma which was obtained by some ideas and methods used in \cite{LYZ}.

\begin{lem}\label{lem:1}
	Suppose that $f\in L^{p}(\mathbb{R}^{n_{1}}\times\mathbb{R}^{n_{2}}\times\cdots\times\mathbb{R}^{n_{m}},|x|^{\gamma})$ with $|x|^{\gamma}=|x_{1}|^{\gamma_{1}}|x_{2}|^{\gamma_{2}}\cdots|x_{m}|^{\gamma_{m}}$. Let
	$$g_{f}(x)=\frac{1}{\omega_{n_{1}}}\frac{1}{\omega_{n_{2}}}\cdots\frac{1}{\omega_{n_{m}}}\int_{|\xi_{1}|=1}\int_{|\xi_{2}|=1}
\cdots\int_{|\xi_{m}|=1}f(|x_{1}|\xi_{1},|x_{2}|\xi_{2},\cdots,|x_{m}|\xi_{m})d\xi_{1}d\xi_{2}\cdots d\xi_{m},$$
	where $x=(x_{1},x_{2},\cdots,x_{m})\in\mathbb{R}^{n_{1}}\times\mathbb{R}^{n_{2}}\times\cdots\times\mathbb{R}^{n_{m}}$, $\omega_{n_{i}}=2\pi^{\frac{n_{i}}{2}}/\Gamma(n_{i}/2)$ with $1\leq i\leq m$. Then
	$$\mathscr{H}_{\beta_{1},\beta_{2},\cdots,\beta_{m}}(|f|)(x)=\mathscr{H}_{\beta_{1},\beta_{2},\cdots,\beta_{m}}(g_{f})(x)$$
	and
	$$\|g_{f}\|_{L^{p}(\mathbb{R}^{n_{1}}\times\mathbb{R}^{n_{2}}\times\cdots\times\mathbb{R}^{n_{m}},|x|^{\gamma})}\leq \|f\|_{L^{p}(\mathbb{R}^{n_{1}}\times\mathbb{R}^{n_{2}}\times\cdots\times\mathbb{R}^{n_{m}},|x|^{\gamma})}.$$
\end{lem}

\noindent$Proof$. We merely the proof with the case $m=2$ for the sake of clarity in writing, and the same is true for the general case $m>2$.

It follows that $\mathscr{H}_{\beta_{1},\beta_{2}}(g_{f})(x_{1},x_{2})$ is equal to
$$
{
	\begin{split}
	&\frac{1}{|B(0,|x_{1}|)|^{1-\frac{\beta_{1}}{n_{1}}}}\frac{1}{|B(0,|x_{2}|)|^{1-\frac{\beta_{2}}{n_{2}}}}\int_{|y_{1}|<|x_{1}|}
\int_{|y_{2}|<|x_{2}|}\frac{1}{\omega_{n_{1}}\omega_{n_{2}}}\\
	&\qquad\qquad\times\int_{|\xi_{1}|=1}\int_{|\xi_{2}|=1}f(|y_{1}|\xi_{1},|y_{2}|\xi_{2})d\xi_{1}d\xi_{2}dy_{1}dy_{2}	\\
	&=\frac{1}{\omega_{n_{1}}\omega_{n_{2}}}\int_{|\xi_{1}|=1}\int_{|\xi_{2}|=1}\frac{1}{|B(0,|x_{1}|)|^{1-\frac{\beta_{1}}{n_{1}}}}
\frac{1}{|B(0,|x_{2}|)|^{1-\frac{\beta_{2}}{n_{2}}}}\\
	&\qquad\qquad\times\int_{|y_{1}|<|x_{1}|}\int_{|y_{2}|<|x_{2}|}f(|y_{1}|\xi_{1},|y_{2}|\xi_{2})dy_{1}dy_{2}d\xi_{1}d\xi_{2}\\
	&=\int_{S^{n_{1}-1}}\int_{S^{n_{2}-1}}\frac{1}{|B(0,|x_{1}|)|^{1-\frac{\beta_{1}}{n_{1}}}}\frac{1}{|B(0,|x_{2}|)|^{1-\frac{\beta_{2}}{n_{2}}}}\\
	&\qquad\times\int_{0}^{|x_{1}|}\int_{0}^{|x_{2}|}f(r_{1}\xi_{1},r_{2}\xi_{2})r_{1}^{n_{1}-1}r_{2}^{n_{2}-1}dr_{1}dr_{2}d\sigma(\xi_{1})
d\sigma(\xi_{2})\\
	&=\mathscr{H}_{\beta_{1},\beta_{2}}(f)(x_{1},x_{2}).
	\end{split}
	}
$$

Using the generalized Minkowski's inequality and H\"{o}lder's inequality, we conclude that

 \noindent$\|g_{f}\|_{L^{p}(\mathbb{R}^{n_{1}}\times\mathbb{R}^{n_{2}},|x|^{\gamma})}$ is not greater than
$$
{
	\begin{split} &\frac{1}{\omega_{n_{1}}\omega_{n_{2}}}\int_{|\xi_{1}|=1}\int_{|\xi_{2}|=1}\left(\int_{\mathbb{R}^{n_{1}}}\int_{\mathbb{R}^{n_{2}}}
\left(f(|x_{1}|\xi_{1},|x_{2}|\xi_{2})\right)^{p}|x_{1}|^{\gamma_{1}}|x_{2}|^{\gamma_{2}}dx_{1}dx_{2}\right)^{\frac{1}{p}}d\xi_{1}d\xi_{2}\\
	&\leq\left(\frac{1}{\omega_{n_{1}}\omega_{n_{2}}}\int_{|\xi_{1}|=1}\int_{|\xi_{2}|=1}\int_{\mathbb{R}^{n_{1}}}\int_{\mathbb{R}^{n_{2}}}
\left(f(|x_{1}|\xi_{1},|x_{2}|\xi_{2})\right)^{p}|x_{1}|^{\gamma_{1}}|x_{2}|^{\gamma_{2}}dx_{1}dx_{2}d\xi_{1}d\xi_{2}\right)^{\frac{1}{p}}\\
	&=\left(\int_{S^{n_{1}-1}}\int_{S^{n_{2}-1}}\int_{0}^{\infty}\int_{0}^{\infty}
\left(f(r_{1}\xi_{1},r_{2}\xi_{2})\right)^{p}
r_{1}^{\gamma_{1}+n_{1}-1}r_{2}^{\gamma_{2}+n_{2}-1}dr_{1}dr_{2}d\sigma(\xi_{1})d\sigma(\xi_{2})\right)^{\frac{1}{p}}\\
	&=\|f\|_{L^{p}(\mathbb{R}^{n_{1}}\times\mathbb{R}^{n_{2}},|x|^{\gamma})}.
	\end{split}
	}
$$
This finishes the proof of the lemma.
$\hfill$ $\Box$

\begin{rem}\label{re:1}
	It follows from Lemma $\ref{lem:1}$ that
	$$\frac{\|\mathscr{H}_{\beta_{1},\beta_{2},\cdots,\beta_{m}}(f)\|_{L^{q}(\mathbb{R}^{n_{1}}\times\mathbb{R}^{n_{2}}
\times\cdots\times\mathbb{R}^{n_{m}},|x|^{\alpha})}}{\|f\|_{L^{p}(\mathbb{R}^{n_{1}}
\times\mathbb{R}^{n_{2}}\times\cdots\times\mathbb{R}^{n_{m}},|x|^{\gamma})}}\leq \frac{\|\mathscr{H}_{\beta_{1},\beta_{2},\cdots,\beta_{m}}(g_{f})\|_{L^{q}(\mathbb{R}^{n_{1}}
\times\mathbb{R}^{n_{2}}\times\cdots\times\mathbb{R}^{n_{m}},|x|^{\alpha})}}{\|g_{f}\|_{L^{p}(\mathbb{R}^{n_{1}}\times\mathbb{R}^{n_{2}}
\times\cdots\times\mathbb{R}^{n_{m}},|x|^{\gamma})}}.$$
	Therefore, the norm of the operator $\mathscr{H}_{\beta_{1},\beta_{2},\cdots,\beta_{m}}$ from $L^{p}(\mathbb{R}^{n_{1}}\times\mathbb{R}^{n_{2}}\times\cdots\times\mathbb{R}^{n_{m}},|x|^{\gamma})$ to $L^{q}(\mathbb{R}^{n_{1}}\times\mathbb{R}^{n_{2}}\times\cdots\times\mathbb{R}^{n_{m}},|x|^{\alpha})$ is equal to the norm that $\mathscr{H}_{\beta_{1},\beta_{2},\cdots,\beta_{m}}$ restricts to radial functions.
\end{rem}

For two differences power weight of fractional Hardy operator $\mathbb{H}_{\beta}$ with $\beta=0$ (write as $\mathbb{H}$), we have following lemma,  which can be found in the paper \cite{PS}.
\begin{lem}\label{lem:2}
	Suppose that $1<p< q<\infty$, $n\in\mathbb{N}$, $x\in\mathbb{R}^{n}$, $\gamma<n(p-1)$ and $\frac{\alpha+n}{q}=\frac{\gamma+n}{p}$. If $f\in L^{p}(\mathbb{R}^{n},|x|^{\gamma})$, then we have
	\begin{equation}
	\|\mathbb{H}(f)\|_{L^{q}(\mathbb{R}^{n},|x|^{\alpha})}\leq C_{pq}^{*}\|f\|_{L^{p}(\mathbb{R}^{n},|x|^{\gamma})},
	\end{equation}
	where $C_{pq}^{*}$ is sharp and equal to
	$$|S^{n-1}|^{\frac{1}{q}-\frac{1}{p}}n^{\frac{1}{p}-\frac{1}{q}}\left(\frac{n(p-1)}{n(p-1)-\gamma}\right)^{\frac{1}{p^{\prime}}+\frac{1}{q}}
\left(\frac{p^{\prime}}{q}\right)^{1/q}\left(\frac{p}{q-p}\cdot B\left(\frac{p}{q-p},\frac{pq}{q^{\prime}(q-p)}\right)\right)^{\frac{1}{q}-\frac{1}{p}}.$$
\end{lem}

\begin{rem}\label{re:2}
	There holds
	$$C_{pq}^{*}\rightarrow\frac{p}{p-1-\frac{\gamma}{n}}$$
	as $q\rightarrow p$.
\end{rem}

In fact, by using Persson and Samko \cite{PS} result:
$$\left(\frac{p^{\prime}}{q}\right)^{1/q}\left(\frac{p}{q-p}\cdot B\left(\frac{p}{q-p},\frac{pq}{q^{\prime}(q-p)}\right)\right)^{\frac{1}{q}-\frac{1}{p}}\rightarrow \frac{p}{p-1}\quad as \quad q\rightarrow p,$$
therefore, we find that
$$C_{pq}^{*}\approx |S^{n-1}|^{\frac{1}{q}-\frac{1}{p}}n^{\frac{1}{p}-\frac{1}{q}}\left(\frac{n(p-1)}{n(p-1)-\gamma}\right)^{\frac{1}{p^{\prime}}+\frac{1}{q}}\frac{p}{p-1}
=\frac{p}{p-1-\frac{\gamma}{n}} \quad as \quad q\rightarrow p.$$

With the help of previous consequences, we shall prove our main statements.
\section{Proof of main results}
\qquad First we use the Lu and Zhao \cite{LZ} result to derive a new constant, which is sharp in $\eqref{eq:4}$ for each $p\in(1,q)$.

\noindent$Proof$ $of$ $Theorem$ $\ref{main:1}$ We merely the proof with the case $m=2$ for the sake of clarity in writing, and the same is true for the general case $m>2$.

Without loss of generality, we suppose that $f$ is an nonnegative  integrable function. For fixed variable $x_{2}$ we denote
$$F(y_{1}):=F(y_{1},x_{2})=\int_{|y_{2}|<|x_{2}|}f(y_{1},y_{2})dy_{2}.$$
According to Lu and Zhao \cite{LZ} estimate for the high dimensional fractional Hardy operator in the case $1<p<q<\infty$ we find that
$$
{
	\begin{split}
	&\int_{\mathbb{R}^{n_{1}}}\int_{\mathbb{R}^{n_{2}}}\left|\frac{1}{|B(0,|x_{1}|)|^{1-\frac{\beta_{1}}{n_{1}}}}
\frac{1}{|B(0,|x_{2}|)|^{1-\frac{\beta_{2}}{n_{2}}}}\int_{|y_{1}|<|x_{1}|}\int_{|y_{2}|<|x_{2}|}f(y_{1},y_{2})dy_{1}dy_{2}\right|^{q}dx_{1}dx_{2}\\
	&=\int_{\mathbb{R}^{n_{2}}}\frac{1}{|B(0,|x_{2}|)|^{\left(1-\frac{\beta_{2}}{n_{2}}\right)q}}\left[\int_{\mathbb{R}^{n_{1}}}
\left(\frac{1}{|B(0,|x_{1}|)|^{1-\frac{\beta_{1}}{n_{1}}}}\int_{|y_{1}|<|x_{1}|}F(y_{1})dy_{1}\right)^{q}dx_{1}\right]dx_{2}\\
	&\leq A_{1}\int_{\mathbb{R}^{n_{2}}}\frac{1}{|B(0,|x_{2}|)|^{\left(1-\frac{\beta_{2}}{n_{2}}\right)q}}
\left(\int_{\mathbb{R}^{n_{1}}}F^{p}(y_{1})dy_{1}\right)^{\frac{q}{p}}dx_{2}=:II,
	\end{split}
	}
$$
where $A_{1}$ is a constant in the inequality of $\eqref{eq:3}$ with $n=n_{1}$.

{\noindent}By applying the generality Minkowski's inequality with the power $\frac{q}{p}$, we obtain that
$$
{
	\begin{split}
	II&\leq A_{1}\left(\int_{\mathbb{R}^{n_{1}}}\left(\int_{\mathbb{R}^{n_{2}}}\frac{1}{|B(0,|x_{2}|)|^{\left(1-\frac{\beta_{2}}{n_{2}}\right)q}}
F^{q}(y_{1})dx_{2}\right)^{\frac{p}{q}}dy_{1}\right)^{\frac{q}{p}}\\
	&=A_{1}\left(\int_{\mathbb{R}^{n_{1}}}\left(\int_{\mathbb{R}^{n_{2}}}\left(\frac{1}{|B(0,|x_{2}|)|^{1-\frac{\beta_{2}}{n_{2}}}}\int_{|y_{2}|<|x_{2}|}
f(y_{1},y_{2})dy_{2}\right)^{q}dx_{2}\right)^{\frac{p}{q}}dy_{1}\right)^{\frac{q}{p}}\\
	&\leq A_{1}A_{2}\left(\int_{\mathbb{R}^{n_{1}}}\int_{\mathbb{R}^{n_{2}}}f^{p}(y_{1},y_{2})dy_{1}dy_{2}\right)^{\frac{q}{p}},
	\end{split}
	}
	$$
where $A_{2}$ is a constant in the inequality of $\eqref{eq:3}$ with $n=n_{2}$.

Therefore, it implies that
\begin{equation}\label{eq:6}
\|\mathscr{H}_{\beta_{1},\beta_{2}}(f)\|_{L^{q}(\mathbb{R}^{n_{1}}\times\mathbb{R}^{n_{2}})}\leq\prod_{i=1}^{2}\left(\frac{p^{\prime}}{q}\right)^{1/q}
\left(\frac{n_{i}}{q\beta_{i}}\cdot B\left(\frac{n_{i}}{q\beta_{i}},\frac{n_{i}}{q^{\prime}\beta_{i}}\right)\right)^{-\beta_{i}/n_{i}}
\|f\|_{L^{p}(\mathbb{R}^{n_{1}}\times\mathbb{R}^{n_{2}})}.
\end{equation}

Next, we need to proved the converse inequality.

It follows from Lemma $\ref{lem:1}$ that the norm of the operator $\mathscr{H}_{\beta_{1},\beta_{2}}$ from $L^{p}(\mathbb{R}^{n_{1}}\times\mathbb{R}^{n_{2}})$ to $L^{q}(\mathbb{R}^{n_{1}}\times\mathbb{R}^{n_{2}})$ is equal to the norm that $\mathscr{H}_{\beta_{1},\beta_{2}}$ restricts to radial functions. Consequently, without loss of generality, it suffices to carry out the proof the converse inequality by assuming that $f$ is a nonnegative, radial, smooth function with compact support on $\mathbb{R}^{n_{1}}\times\mathbb{R}^{n_{2}}$.

Using the polar coordinate transformation, we can rewrite $\eqref{eq:6}$ as
\begin{equation}\label{eq:8}
{
	\begin{split}
	&n_{1}\int_{0}^{\infty}n_{2}\int_{0}^{\infty}
\left(n_{1}\int_{0}^{s_{1}}n_{2}\int_{0}^{s_{2}}f(r_{1},r_{2})r_{1}^{n_{1}-1}r_{2}^{n_{2}-1}dr_{1}dr_{2}\right)^{q}
s_{1}^{q(\beta_{1}-n_{1})}s_{1}^{n_{1}-1}s_{2}^{q(\beta_{2}-n_{2})}s_{2}^{n_{2}-1}ds_{1}ds_{2}\\
	&\leq\left(A_{1}A_{2}\right)^{q}\left(n_{1}\int_{0}^{\infty}n_{2}\int_{0}^{\infty}f^{p}
(s_{1},s_{2})s_{1}^{n_{1}-1}s_{2}^{n_{2}-1}ds_{1}ds_{2}\right)^{\frac{q}{p}}.
	\end{split}
	}
\end{equation}

For the purpose of getting the sharp bound, we take $\tilde{f}(x_{1},x_{2})=\frac{1}{\left(1+|x_{1}|^{q\beta_{1}}\right)^{1+\frac{n_{1}}{q\beta_{1}}}}
\frac{1}{\left(1+|x_{2}|^{q\beta_{2}}\right)^{1+\frac{n_{2}}{q\beta_{2}}}}$. It follows from
$$n_{1}\int_{0}^{s_{1}}n_{2}\int_{0}^{s_{2}}\tilde{f}(r_{1},r_{2})r_{1}^{n_{1}-1}r_{2}^{n_{2}-1}dr_{1}dr_{2}=\frac{s_{1}^{n_{1}}}
{\left(1+s_{1}^{q\beta_{1}}\right)^{\frac{n_{1}}{q\beta_{1}}}}\frac{s_{2}^{n_{2}}}{\left(1+s_{2}^{q\beta_{2}}\right)^{\frac{n_{2}}{q\beta_{2}}}}$$
that the left side of $\eqref{eq:8}$ is
$$\prod_{i=1}^{2}\frac{n_{i}}{q\beta_{i}}\cdot B\left(\frac{n_{i}}{q\beta_{i}}+1,\frac{n_{i}}{q^{\prime}\beta_{i}}-1\right).$$
It is easy to verify that
$$n_{1}\int_{0}^{\infty}n_{2}\int_{0}^{\infty}\tilde{f}^{p}(s_{1},s_{2})s_{1}^{n_{1}-1}s_{2}^{n_{2}-1}ds_{1}ds_{2}
=\prod_{i=1}^{2}\frac{n_{i}}{q\beta_{i}}\cdot B\left(\frac{n_{i}}{q\beta_{i}},\frac{n_{i}}{q^{\prime}\beta_{i}}\right).$$
Therefore,
$$
{
	\begin{split}
	\|\mathscr{H}_{\beta_{1},\beta_{2}}\|_{L^{p}(\mathbb{R}^{n_{1}}\times\mathbb{R}^{n_{2}})\rightarrow L^{q}(\mathbb{R}^{n_{1}}\times\mathbb{R}^{n_{2}})}&=\sup_{\|f\|_{L^{p}(\mathbb{R}^{n_{1}}\times\mathbb{R}^{n_{2}})}\neq0}
\frac{\|\mathscr{H}_{\beta_{1},\beta_{2}}(f)\|_{L^{q}(\mathbb{R}^{n_{1}}\times\mathbb{R}^{n_{2}})}}{\|f\|_{L^{p}(\mathbb{R}^{n_{1}}
\times\mathbb{R}^{n_{2}})}}\\
	&\geq\frac{\|\mathscr{H}_{\beta_{1},\beta_{2}}(\tilde{f})\|_{L^{q}(\mathbb{R}^{n_{1}}\times\mathbb{R}^{n_{2}})}}{\|\tilde{f}\|_{L^{p}
(\mathbb{R}^{n_{1}}\times\mathbb{R}^{n_{2}})}}=A_{1}A_{2},
	\end{split}
	}
$$
This completes the proof of Theorem $\ref{main:1}$.   $\hfill$ $\Box$

\noindent$Proof$ $of$ $Theorem$ $\ref{main:2}$  We merely the proof with the case $m=2$ for the sake of clarity in writing, and the same is true for the general case $m>2$.

Without loss of generality, it follows from Lemma $\ref{lem:1}$ we can assuming that $f$ is a nonnegative, radial, smooth function with compact support on $\mathbb{R}^{n_{1}}\times\mathbb{R}^{n_{2}}$.

Using the polar coordinate transformation, $\|\mathscr{H}_{\beta_{1},\beta_{2}}(f)\|_{L^{q}(\mathbb{R}^{n_{1}}\times\mathbb{R}^{n_{2}},|x|^{\alpha})}$ is equal to
\begin{equation}\label{eq:9}
{
	\begin{split}
	&K_{q}\bigg(\int_{0}^{\infty}\int_{0}^{\infty}\left(\int_{0}^{\rho_{1}}\int_{0}^{\rho_{2}}
f(t_{1},t_{2})t_{1}^{n_{1}-1}t_{2}^{n_{2}-1}dt_{1}dt_{2}\right)^{q}\\
	&\quad\qquad\qquad\times\rho_{1}^{q(\beta_{1}-n_{1})+\alpha_{1}}\rho_{1}^{n_{1}-1}
\rho_{2}^{q(\beta_{2}-n_{2})+\alpha_{2}}\rho_{2}^{n_{2}-1}d\rho_{1}d\rho_{2}\bigg)^{\frac{1}{q}}\\
	\end{split}
}
\end{equation}
and $\|f\|_{L^{p}(\mathbb{R}^{n_{1}}\times\mathbb{R}^{n_{2}},|x|^{\gamma})}$ is equal to
\begin{equation}\label{eq:10}
K_{p}\left(\int_{0}^{\infty}\int_{0}^{\infty}f^{p}(\rho_{1},\rho_{2})
\rho_{1}^{n_{1}-1+\gamma_{1}}\rho_{2}^{n_{2}-1+\gamma_{2}}d\rho_{1}d\rho_{2}\right)^{\frac{1}{p}},
\end{equation}
where two constants $K_{q}=n_{1}^{1-\frac{\beta_{1}}{n_{1}}}n_{2}^{1-\frac{\beta_{2}}{n_{2}}}|S^{n_{1}-1}|^{\frac{1}{q}
+\frac{\beta_{1}}{n_{1}}}|S^{n_{2}-1}|^{\frac{1}{q}+\frac{\beta_{2}}{n_{2}}}$ and $K_{p}=|S^{n_{1}-1}|^{\frac{1}{p}}|S^{n_{2}-1}|^{\frac{1}{p}}$.

First we make a change of variables in $\eqref{eq:10}$ by putting
\begin{equation}\label{eq:11}
\frac{s_{1}^{n_{1}}}{n_{1}}=\frac{s_{1}^{n_{1}}(\rho_{1})}{n_{1}}=\int_{0}^{\rho_{1}}t_{1}^{n_{1}-1
-\frac{\gamma_{1}}{p-1}}dt=\frac{p-1}{n_{1}(p-1)-\gamma_{1}}\rho_{1}^{\frac{n_{1}(p-1)-\gamma_{1}}{p-1}},
\end{equation}

\begin{equation}\label{eq:12}
\frac{s_{2}^{n_{2}}}{n_{2}}=\frac{s_{2}^{n_{2}}(\rho_{2})}{n_{2}}=\int_{0}^{\rho_{2}}t_{2}^{n_{2}-1
-\frac{\gamma_{2}}{p-1}}dt=\frac{p-1}{n_{2}(p-1)-\gamma_{2}}\rho_{2}^{\frac{n_{2}(p-1)-\gamma_{2}}{p-1}}
\end{equation}
and define
\begin{equation}\label{eq:13}
g(s_{1},s_{2})=g(s_{1}(\rho_{1}),s_{2}(\rho_{2}))=f(\rho_{1},\rho_{2})\rho_{1}^{\frac{\gamma_{1}}{p-1}}\rho_{2}^{\frac{\gamma_{2}}{p-1}}.
\end{equation}
Then
\begin{equation}\label{eq:14}
{
	\begin{split}
	&\int_{0}^{\infty}\int_{0}^{\infty}f^{p}(\rho_{1},\rho_{2})\rho_{1}^{n_{1}-1+\gamma_{1}}\rho_{2}^{n_{2}-1+\gamma_{2}}d\rho_{1}d\rho_{2}\\
	&=\int_{0}^{\infty}\int_{0}^{\infty}f^{p}(\rho_{1},\rho_{2})\rho_{1}^{\frac{\gamma_{1}p}{p-1}}\rho_{1}^{n_{1}-1-\frac{\gamma_{1}}{p-1}}d\rho_{1}
\rho_{2}^{\frac{\gamma_{2}p}{p-1}}\rho_{2}^{n_{2}-1-\frac{\gamma_{2}}{p-1}}d\rho_{2}\\
	&=\int_{0}^{\infty}\int_{0}^{\infty}g^{p}(s_{1},s_{2})d\frac{s_{1}^{n_{1}}}{n_{1}}d\frac{s_{2}^{n_{2}}}{n_{2}}
=\int_{0}^{\infty}\int_{0}^{\infty}g^{p}(s_{1},s_{2})s_{1}^{n_{1}-1}s_{2}^{n_{2}-1}ds_{1}ds_{2}
	\end{split}
	}
\end{equation}
and
\begin{equation}\label{eq:15}
{
	\begin{split}
	II_{0}:=&\bigg(\int_{0}^{\infty}\int_{0}^{\infty}\left(\int_{0}^{\rho_{1}}
\int_{0}^{\rho_{2}}f(t_{1},t_{2})t_{1}^{n_{1}-1}t_{2}^{n_{2}-1}dt_{1}dt_{2}\right)^{q}\\
	&\quad\qquad\qquad\times\rho_{1}^{q(\beta_{1}-n_{1})+\alpha_{1}}\rho_{1}^{n_{1}-1}\rho_{2}^{q(\beta_{2}-n_{2})
+\alpha_{2}}\rho_{2}^{n_{2}-1}d\rho_{1}d\rho_{2}\bigg)^{\frac{1}{q}}\\
	&=\bigg(\int_{0}^{\infty}\int_{0}^{\infty}\left(\int_{0}^{\rho_{1}}\int_{0}^{\rho_{2}}
f(t_{1},t_{2})t_{1}^{n_{1}-1}t_{2}^{n_{2}-1}dt_{1}dt_{2}\right)^{q}\\
	&\quad\times\rho_{1}^{q(\beta_{1}-n_{1})+\left(\frac{\gamma_{1}+n_{1}}{p}-\beta_{1}\right)q-n_{1}}\rho_{1}^{n_{1}-1}\rho_{2}^{q(\beta_{2}-n_{2})
+\left(\frac{\gamma_{2}+n_{2}}{p}-\beta_{2}\right)q-n_{2}}\rho_{2}^{n_{2}-1}d\rho_{1}d\rho_{2}\bigg)^{\frac{1}{q}}.
	\end{split}
	}
\end{equation}
To obtain we desire result, we need to following form. And $\eqref{eq:15}$ is equal to
\begin{equation*}
{
	\begin{split}
	II_{0}&=\bigg(\int_{0}^{\infty}\int_{0}^{\infty}\left(\int_{0}^{\rho_{1}}\int_{0}^{\rho_{2}}f(t_{1},t_{2})
t_{1}^{n_{1}-1}t_{2}^{n_{2}-1}dt_{1}dt_{2}\right)^{q}\rho_{1}^{q(\beta_{1}-n_{1})+
\left(\frac{\gamma_{1}+n_{1}}{p}-\beta_{1}\right)q-n_{1}+\frac{\gamma_{1}}{p-1}}\\
	&\qquad\qquad\qquad\times\rho_{1}^{n_{1}-1-\frac{\gamma_{1}}{p-1}}\rho_{2}^{q(\beta_{2}-n_{2})+
\left(\frac{\gamma_{2}+n_{2}}{p}-\beta_{2}\right)q-n_{2}+\frac{\gamma_{2}}{p-1}}\rho_{2}^{n_{2}-1-
\frac{\gamma_{2}}{p-1}}d\rho_{1}d\rho_{2}\bigg)^{\frac{1}{q}}.
	\end{split}
	}
\end{equation*}
Hence, since
\begin{equation}\label{eq:16}
{\begin{split}
f(t_{1},t_{2})t_{1}^{n_{1}-1}t_{2}^{n_{2}-1}dt_{1}dt_{2}&=f(t_{1},t_{2})t_{1}^{\frac{\gamma_{1}}{p-1}}
t_{1}^{n_{1}-1-\frac{\gamma_{1}}{p-1}}dt_{1}t_{2}^{\frac{\gamma_{2}}{p-1}}t_{2}^{n_{2}-1-\frac{\gamma_{2}}{p-1}}dt_{2}\\
&=g(r_{1},r_{2})r_{1}^{n_{1}-1}r_{2}^{n_{2}-1}dr_{1}dr_{2}
\end{split}
}
\end{equation}
and
\begin{equation}\label{eq:17}
\frac{s_{1}^{n_{1}}}{n_{1}}=\frac{p-1}{n_{1}(p-1)-\gamma_{1}}\rho_{1}^{\frac{n_{1}(p-1)-\gamma_{1}}{p-1}}\Rightarrow \rho_{1}=\left(\frac{n_{1}(p-1)-\gamma_{1}}{n_{1}(p-1)s_{1}^{-n_{1}}}\right)^{\frac{p-1}{n_{1}(p-1)-\gamma_{1}}},
\end{equation}
similarly,
\begin{equation}\label{eq:18}
\rho_{2}=\left(\frac{n_{2}(p-1)-\gamma_{2}}{n_{2}(p-1)s_{2}^{-n_{2}}}\right)^{\frac{p-1}{n_{2}(p-1)-\gamma_{2}}}.
\end{equation}
Applying the $\eqref{eq:15}$, $\eqref{eq:16}$, $\eqref{eq:17}$, and $\eqref{eq:18}$ we have that
\begin{equation}\label{eq:19}
{
	\begin{split}
	II_{0}=&\prod_{i=1}^{2}\left(\frac{n_{i}(p-1)}{n_{i}(p-1)-\gamma_{i}}\right)^{\frac{1}{p^{\prime}}+\frac{1}{q}}
\bigg(\int_{0}^{\infty}\int_{0}^{\infty}\left(\int_{0}^{s_{1}}\int_{0}^{s_{2}}g(r_{1},r_{2})r_{1}^{n_{1}-1}r_{2}^{n_{2}-1}dr_{1}dr_{2}\right)^{q}\\
	&\qquad\qquad\qquad\qquad\qquad\qquad\times s_{1}^{-n_{1}(\frac{q}{p^{\prime}}+1)}s_{1}^{n_{1}-1}s_{2}^{-n_{2}(\frac{q}{p^{\prime}}+1)}s_{2}^{n_{2}-1}ds_{1}ds_{2}\bigg)^{\frac{1}{q}}.
	\end{split}
	}
\end{equation}
Taking $\beta_{1}=n_{1}(\frac{1}{p}-\frac{1}{q})$ and $\beta_{2}=n_{2}(\frac{1}{p}-\frac{1}{q})$ in $\eqref{eq:8}$, we obtain that

\begin{equation}\label{eq:20}
{
	\begin{split}
	&n_{1}\int_{0}^{\infty}n_{2}\int_{0}^{\infty}
\left(n_{1}\int_{0}^{s_{1}}n_{2}\int_{0}^{s_{2}}f(r_{1},r_{2})r_{1}^{n_{1}-1}r_{2}^{n_{2}-1}dr_{1}dr_{2}\right)^{q}s_{1}^{-
\frac{n_{1}q}{p^{\prime}}-1}s_{2}^{-\frac{n_{2}q}{p^{\prime}}-1}ds_{1}ds_{2}\\
	&\leq\left(A_{1}A_{2}\right)^{q}\left(n_{1}\int_{0}^{\infty}n_{2}\int_{0}^{\infty}f^{p}
(s_{1},s_{2})s_{1}^{n_{1}-1}s_{2}^{n_{2}-1}ds_{1}ds_{2}\right)^{\frac{q}{p}}.
	\end{split}
}
\end{equation}
Combining $\eqref{eq:19}$ and $\eqref{eq:20}$ we find that
\begin{equation}\label{eq:21}
{
	\begin{split}
	II_{0}\leq \prod_{i=1}^{2}\left(\frac{n_{i}(p-1)}{n_{i}(p-1)
-\gamma_{i}}\right)^{\frac{1}{p^{\prime}}+\frac{1}{q}}n_{i}^{\frac{1}{p}-\frac{1}{q}-1}
A_{i}\left(\int_{0}^{\infty}\int_{0}^{\infty}g^{p}(s_{1},s_{2})s_{1}^{n_{1}-1}s_{2}^{n_{2}-1}ds_{1}ds_{2}\right)^{\frac{1}{p}}.
	\end{split}
	}
\end{equation}
So using the Lemma $\ref{lem:1}$ we implies that
\begin{equation}\label{eq:22}
{
	\begin{split}
	&\frac{\|\mathscr{H}_{\beta_{1},\beta_{2}}(f)\|_{L^{q}(\mathbb{R}^{n_{1}}\times\mathbb{R}^{n_{2}},|x|^{\alpha})}}
{\|f\|_{L^{p}(\mathbb{R}^{n_{1}}\times\mathbb{R}^{n_{2}},|x|^{\gamma})}}\leq \frac{\|\mathscr{H}_{\beta_{1},\beta_{2}}(g_{f})\|_{L^{q}(\mathbb{R}^{n_{1}}\times\mathbb{R}^{n_{2}},|x|^{\alpha})}}{\|g_{f}\|_{L^{p}
(\mathbb{R}^{n_{1}}\times\mathbb{R}^{n_{2}},|x|^{\gamma})}}\\
	&\qquad\quad=\frac{K_{q}II_{0}}{K_{p}\left(\int_{0}^{\infty}\int_{0}^{\infty}g^{p}(s_{1},s_{2})s_{1}^{n_{1}-1}s_{2}^{n_{2}-1}
ds_{1}ds_{2}\right)^{\frac{1}{p}}}.
	\end{split}
	}
\end{equation}
Therefore, combining $\eqref{eq:21}$ and $\eqref{eq:22}$, we have $\|\mathscr{H}_{\beta_{1},\beta_{2}}\|_{L^{p}(\mathbb{R}^{n_{1}}\times\mathbb{R}^{n_{2}},|x|^{\gamma})\rightarrow L^{q}(\mathbb{R}^{n_{1}}\times\mathbb{R}^{n_{2}},|x|^{\alpha})}$ is not greater than
\begin{equation*}
\prod_{i=1}^{2}|S^{n_{i}-1}|^{\frac{1}{q}-\frac{1}{p}+\frac{\beta_{i}}{n_{i}}}n_{i}^{\frac{1}{p}-\frac{1}{q}-\frac{\beta_{i}}{n_{i}}}
\left(\frac{n_{i}(p-1)}{n_{i}(p-1)-\gamma_{i}}\right)^{\frac{1}{p^{\prime}}+\frac{1}{q}}\left(\frac{p^{\prime}}{q}\right)^{1/q}\left(\frac{p}{q-p}\cdot B\left(\frac{p}{q-p},\frac{pq}{q^{\prime}(q-p)}\right)\right)^{\frac{1}{q}-\frac{1}{p}}.
\end{equation*}
From the Theorem $\ref{main:1}$ it also follows that above constant is sharp when
$$g(s_{1},s_{2})=\frac{1}{\left(1+s_{1}^{\frac{n_{1}q}{p}-n_{1}}\right)^{\frac{q}{q-p}}}\frac{1}{\left(1+s_{2}^{\frac{n_{2}q}{p}-n_{2}}
\right)^{\frac{q}{q-p}}}$$
i.e. when (see $\eqref{eq:11}$, $\eqref{eq:12}$ and $\eqref{eq:13}$)
$$f(x_{1},x_{2})=\frac{|x_{1}|^{-\frac{\gamma_{1}}{p-1}}}{\left(1+d_{1}|x_{1}|^{\frac{n_{1}(p-1)
-\gamma_{1}}{n_{1}(p-1)}(\frac{n_{1}q}{p}-n_{1})}\right)^{\frac{q}{q-p}}}\frac{|x_{2}|^{-\frac{\gamma_{2}}{p-1}}}
{\left(1+d_{2}|x_{2}|^{\frac{n_{2}(p-1)-\gamma_{2}}{n_{2}(p-1)}(\frac{n_{2}q}{p}-n_{2})}\right)^{\frac{q}{q-p}}},$$
where two constants $d_{1}$ and $d_{2}$ are $\left(\frac{n_{1}(p-1)}{n_{1}(p-1)-\gamma_{1}}\right)^{1/n_{1}}$ and $\left(\frac{n_{2}(p-1)}{n_{2}(p-1)-\gamma_{2}}\right)^{1/n_{2}}$, respectively.

This completes the proof of Theorem $\ref{main:2}$.  $\hfill$ $\Box$

\begin{rem}
	If $\beta_{1}=\beta_{2}=\cdots=\beta_{m}=0$, the constant is also sharp  in Theorem $\ref{main:2}$. Moreover, the constant is equal to
	$$\prod_{i=1}^{m}|S^{n_{i}-1}|^{\frac{1}{q}-\frac{1}{p}}n_{i}^{\frac{1}{p}-\frac{1}{q}}\left(\frac{n_{i}(p-1)}{n_{i}(p-1)
-\gamma_{i}}\right)^{\frac{1}{p^{\prime}}+\frac{1}{q}}\left(\frac{p^{\prime}}{q}\right)^{1/q}\left(\frac{p}{q-p}\cdot B\left(\frac{p}{q-p},\frac{pq}{q^{\prime}(q-p)}\right)\right)^{\frac{1}{q}-\frac{1}{p}},$$
and it shows that the constant in $\eqref{different power}$ is not really the best.

{\noindent}Using the Lemma $\ref{lem:2}$ and Remark $\ref{re:2}$, we find that
	$$
	{
		\begin{split}
		&\prod_{i=1}^{m}|S^{n_{i}-1}|^{\frac{1}{q}-\frac{1}{p}}n_{i}^{\frac{1}{p}-\frac{1}{q}}\left(\frac{n_{i}(p-1)}{n_{i}(p-1)-\gamma_{i}}
\right)^{\frac{1}{p^{\prime}}+\frac{1}{q}}\left(\frac{p^{\prime}}{q}\right)^{1/q}\left(\frac{p}{q-p}\cdot B\left(\frac{p}{q-p},\frac{pq}{q^{\prime}(q-p)}\right)\right)^{\frac{1}{q}-\frac{1}{p}}\\
		&\qquad\qquad\qquad\qquad\qquad \qquad\approx \prod_{i=1}^{m}\frac{p}{p-1-\alpha_{i}/n_{i}}\quad as \quad p\rightarrow q
		\end{split}
		}
	$$
	and it remains to note that $\prod_{i=1}^{m}\frac{p}{p-1-\alpha_{i}/n_{i}}$ is the sharp constant in \rm\cite{LYZ}.
\end{rem}

{\noindent}\textrm{\bf Acknowledgements} \quad This work was supported by National Natural Science Foundation of China (Grant Nos. 11471309 and 11561062). The authors cordially thank the referees and the editors for their careful reading of the paper and valuable comments which led to the great improvement of the paper.

\end{document}